\documentclass[12pt]{amsart}
\usepackage[centertags]{amsmath}

\usepackage{amsfonts}
\usepackage{amssymb}

\title{A remark on a theorem of Erd\H{o}s}

\author{james H.~Schmerl}
\date{\today}

\def\into{\longrightarrow}

\DeclareMathOperator{\dom}{dom}

\DeclareMathOperator{\supp}{supp}
\DeclareMathOperator{\cf}{cf}

\def\.{.\@}

\def\res{\hspace{-4pt} \upharpoonright \hspace{-2pt}}
\def\harp{\res}

\def\RR{\mathbb{R}} 

\def\FF{\mathbb{F}}

\def\into{\longrightarrow}

\begin{document}
\vspace{-1in}

\begin{abstract} A theorem of Erd\H{o}s  asserts that  every infinite $X \subseteq \RR^n$ has a subset of the same cardinality having no repeated distances. This theorem is generalized here as follows: If $(\RR^n,E)$ is an algebraic hypergraph that does not have an infinite, complete subset,   then every infinite subset has an independent subset of the same cardinality.
\end{abstract}

\maketitle

Erd\H{o}s published the following theorem twice. It first appeared in 1950 in \cite[Th.~III]{e50} 
and then reappeared  in 1978/79 in \cite[Theorem~1]{e} where he remarks  that he had 
proved  it about 
30 years earlier, but ``the published  proof is obscure and not accurate.''

\bigskip

{\sc Theorem 1}: {\em Suppose that $1 \leq n < \omega$. Every infinite  $X \subseteq \RR^n$ has a subset $Y \subseteq X$ 
such that $|Y| = |X|$ and $Y$ has no repeated nonzero distances $($i.e. whenever   $a,b,c,d \in Y$ are such that 
$\|a-b\| = \|c-d\| > 0$, then  $\{a,b\} = \{c,d\})$.} 

\bigskip

 Given a ($k \times n$)-ary polynomial{\footnote{All polynomials considered here are tacitly understood to be over the  reals $\RR$.}}
 $p(x_0,x_1, \ldots,x_{k-1})$ (so that each $x_i$ is an $n$-tuple of variables) and $X \subseteq \RR^n$,  say that 
 $X$ is {\bf independent} (for $p(x_0,x_1, \ldots,x_{k-1})$) if whenever $a_0,a_1, \ldots,$ $a_{k-1} \in X$ are distinct, then $p(a_0,a_1, \ldots, a_{k-1}) \neq 0$.
Let 
$p_n(x,y,z,w)$ be the $(4\times n)$-ary polynomial $\|x - y\|^2 - \|z - w\|^2$, and let 
$q_n(x,y,z)$ be the $(3\times n)$-ary polynomial $p_n(x,y,z,y)$. For any infinite $X \subseteq \RR^n$, $X$ has no repeated distances iff $X$ is independent for both $p_n(x,y,z,w)$ and $q_n(x,y,z)$ or, equivalently, for their product $p_n(x,y,z,w) \cdot q_n(x,y,z)$. Hence, Theorem~1 can be rephrased as follows:
 {\em  Every infinite $X \subseteq \RR^n$ has a subset $Y \subseteq X$ 
such that $|Y| = |X|$ and $Y$ is independent for  $p_n(x,y,z,w) \cdot q_n(x,y,z)$.}  This suggests the question: Which   $(k \times n)$-ary polynomials $p(x_0,x_1, \ldots, x_{k-1})$ have the property that every infinite $X \subseteq \RR^n$ has an independent  subset $Y \subseteq X$ such that $|Y| = |X|$? That question is answered  here in Theorem~2.

It is perhaps more apt to recast this question  in terms of hypergraphs. For any set $X$ and $k < \omega$, we let $[X]^k$ be the set of $k$-element subsets of~$X$. If $2 \leq k < \omega$, then a $k$-{\bf hypergraph} is a pair $H = (V,E)$, where $V$ is its set of vertices and $E \subseteq [V]^k$ its  set of edges.
By a hypergraph, we mean a $k$-hypergraph, where $2 \leq k < \omega$.{\footnote{ What we are calling $k$-hypergraphs are often referred to as $k$-uniform hypergraphs. Thus, any  hypergraph referred to here is 
understood to be a uniform hypergraph.}}
If $H = (V,E)$ is a $k$-hypergraph and $X \subseteq V$, then $X$ is {\bf independent} if $[X]^k \cap E = \varnothing$, and $X$ is {\bf complete} if $[X]^k \subseteq E$.   A $k$-hypergraph $H = (V,E)$ is {\bf algebraic} if there are $n < \omega$ and a $(k \times n)$-ary polynomial $p(x_0,x_1, \ldots,x_{k-1})$ such that 
$V = \RR^n$ and $E$ is the set of those  $\{a_0,a_1, \ldots, a_{k-1}\} \in [\RR^n]^k$  such that 
$p(a_0,a_1, \ldots, a_{k-1}) = 0$.  We say that such $p(x_0,x_1, \ldots, x_{k-1})$ {\bf defines} $H$.

Thus,  Theorem~1 asserts that certain algebraic $4$-hypergraphs, namely those defined by  
$p_n(x,y,z,w)  \cdot q_n(x,y,z)$,   are such that every infinite  set of vertices   has an independent subset  of the same cardinality. What other algebraic hypergraphs have this property? Certainly, any hypergraph having an infinite 
complete set of vertices would not. It is easily seen (by Ramsey's Theorem for $4$-element sets) that none of the hypergraphs defined by 
$p_n(x,y,z,w) \cdot q_n(x,y,z)$ has infinite (or even arbitrarily large finite) complete sets of vertices. Thus, the following theorem, our principal result, generalizes Theorem~1 by showing that the existence of an infinite complete set is the only obstacle.

 \bigskip

{\sc Theorem 2}: {\em Suppose that  $ n < \omega$ and $H = (\RR^n,E)$ is an algebraic hypergraph.  The following are equivalent$:$

\begin{itemize}

\item [(1)] There is a finite $m$ such that whenever  $X \subseteq \RR^n$ is complete, then $|X|<m$.

\item[(2)] If $X \subseteq \RR^n$ is complete, then  $|X| < 2^{\aleph_0}$.

\item[(3)] Every infinite $X \subseteq \RR^n$ has an independent subset $Y \subseteq X$ 
such that $|Y| = |X|$.

\end{itemize} }

\bigskip

The implications $(3) \Longrightarrow  (2) \Longrightarrow (1)$ are known; 
the new part of Theorem~2 is the implication $(1) \Longrightarrow (3)$. The proof of this implication relies on results from \cite{sch16}. We take this opportunity to prove some additional results about algebraic hypergraphs that follow without too much trouble from what is in \cite{sch16}. In \S1 we make explicit  what we are calling a {\it deconstruction} of algebraic hypergraphs that was only implicit in  \cite{sch16}. In \S2 we review some material from \cite{sch16} about chromatic numbers of algebraic hypergraphs, and then make some new observations about the chromatic numbers of their induced subhypergraphs and also about getting a single proper coloring for many hypergraphs at once.
The proof of Theorem~2 is then completed in \S3. 

\bigskip

%
%

{\bf \S1.~The deconstruction of algebraic hypergraphs.} In this section we will make explicit in Theorem~1.1 a result that was previously only implicit in  \cite{sch16}. 

Recall that if $H_0 = (V_0,E_0)$ and $H_1 = (V_1, E_1)$ are two $k$-hypergraphs, then $H_0$ is a {\bf subhypergraph} of $H_1$ if $H_0 \subseteq H_1$ and $E_0 \subseteq E_1$. If, in addition, $V_0 = V_1$, then $H_0$ is a {\bf spanning} subhypergraph of $H_1$, and if $E_0 = E_1 \cap [V_0]^k$, then $H_0$ is an {\bf induced} subhypergraph of $H_1$, in which case we write $H_0 = H_1[V_0]$.  

We will define templates and their hypergraphs as was previously done in \cite[\S2]{sch16}. 
  If $d < \omega$ and $2 \leq k < \omega$, then a $d$-{\bf dimensional} $k$-{\bf template} is a set $P$ of  $d$-tuples such that $|P| = k$. If $P,Q$ are $d$-dimensional $k$-templates, then $f : P \into Q$ is a {\bf homomorphism} from $P$ to $Q$ if $f$ is one-to-one and whenever $x,y \in P$, $i < d$ and $x_i = y_i$, then $f(x)_i = f(y)_i$. 
  If there is such a homomorphism, then $Q$ is a {\bf homomorphic image} of $P$. 
  The  $d$-dimensional $k$-templates $P, Q$ are {\bf isomorphic} if there 
  is a one-to-one function $f : P \into Q$ such that whenever $x,y \in P$ and $i < d$, then    $f(x)_i = f(y)_i$ 
  iff $x_i = y_i$.

  Suppose that $P$ is a $d$-dimensional $k$-template. Let $X_0,X_1, \ldots, X_{d-1}$ be arbitrary sets, and let $X = X_0 \times X_1 \times \cdots \times X_{d-1}$.  Then 
  the {\bf template} hypergraph $L(X,P)$  is the $k$-hypergraph $(X,F)$, where $F$ is the set 
   of those $Q \in [X]^k$ that are homomorphic images of $P$.

   The two theorems of this section make use of  a mild generalization of template hypergraphs. 
   Suppose that $X = X_0 \times X_1 \times \cdots \times X_{d-1}$ and that ${\mathcal P}$ is a {\em set} of $d$-dimensional $k$-templates. Then the {\bf multi-template} hypergraph $L(X,{\mathcal P})$ is the $k$-hypergraph $(X,F)$, where $F$ is the set of 
   those $Q \in [X]^k$ that are homomorphic images of some $P \in {\mathcal P}$. In other words, $H$ is a multi-template hypergraph iff it is the union of finitely many spanning subhypergraphs each of which is a template hypergraph.     Notice that an empty ${\mathcal P}$ is allowed. 
   All multi-template hypergraphs with vertex-set $\RR^d$ are algebraic, and, as such, are especially simple ones. The next theorem shows that they are also quite fundamental. 
     
    \bigskip
   
   {\sc Theorem~1.1}: {\em Suppose that $H= (\RR^n, E)$ is an algebraic $k$-hypergraph. There is a countable collection ${\mathcal A}$ of sets  such that $\RR^n = \bigcup {\mathcal A}$ and for  each $A \in {\mathcal A}$, $H[A]$ is isomorphic to an  algebraic multi-template $k$-hypergraph.}
   
   \bigskip
   
   The previous theorem shows that every algebraic $k$-hypergraph can be deconstructed into countably many  $k$-hypergraphs that are isomorphic to especially simple algebraic $k$-hypergraphs. 
   The next theorem says that Theorem~1.1 can be generalized by considering many algebraic hypergraphs simultaneously, provided we allow a commensurately larger number of subsets.   Theorem~1.1 is a consequence of Theorem~1.2, so we will forgo proving Theorem~1.1 and then  give a sketch of the proof of Theorem~1.2. 
   
   \bigskip
   
   {\sc Theorem 1.2}: {\em Suppose that $1 \leq n < \omega$ and $\aleph_0 \leq \kappa \leq 2^{\aleph_0}$. Suppose that ${\mathcal H}$ is a set of algebraic hypergraphs $H = (\RR^n ,E)$ such that $|{\mathcal H}| \leq \kappa$. Then there  is a collection 
   ${\mathcal A}$ of sets such that $|{\mathcal A}| \leq \kappa$, $\RR^n = \bigcup {\mathcal A}$ and 
   for each $A \in {\mathcal A}$, there are $d < \omega$ and $g : A \into \RR^d$ such that for all $H \in {\mathcal H}$, $g$ is an isomorphism from $H[A]$ to some  algebraic multi-template $k$-hypergraph $L(\RR^d, {\mathcal P})$.}
   
   \bigskip
   
   {\it Proof}.     (We closely follow the proof of $(3) \Longrightarrow (1)$ of Theorem~3.2 in \cite{sch16}.) Let $n, \kappa$ and ${\mathcal H}$ be as given. If $\kappa = 2^{\aleph_0}$, then the existence of ${\mathcal A}$ is trivial as we can let it be the set of all  singleton subsets of $\RR^n$. So assume that $\kappa < 2^{\aleph_0}$. Let $\FF \subseteq \RR$ be a real-closed subfield of $\RR$ such that $|\FF| = \kappa$ and each hypergraph in ${\mathcal H}$ is defined by a polynomial all of whose coefficients are in $\FF$.
   We can now assume that ${\mathcal H}$ is just the set of all algebraic hypergraphs $(\RR^n, E)$ defined by polynomials over $\FF$.

   Let $T$ be a transcendence basis for $\RR$ over $\FF$ such that whenever $x <y \in \RR$, then 
   $|(x,y) \cap T| = 2^{\aleph_0}$. For $a \in \RR^n$, let $\supp(a)$, the {\bf support} of $a$,  be the smallest subset $S \subseteq T$ such that $a$ is algebraic over $\FF \cup S$. 
   
   Suppose $a \in \RR^n$. We say that $f$ is a {\bf determining} function for $a$  if $\supp(a) = \{t_0,t_1, \ldots, t_{d-1}\}$  and    there are 
    rationals $q_0, q_1, \ldots, q_{d-1}, r_0,r_1,$ $ \ldots, r_{d-1}$
  such that:
  
  \begin{itemize}
  
  \item 
   $q_0 < t_0 < r_0 < q_1 < t_1 < \cdots <  q_{d-1} < t_{d-1} < r_{d-1}$;
  
  \item $f : (q_0, r_0) \times (q_1, r_1) \times \cdots \times   (q_{d-1}, r_{d-1}) \into \RR^n$;
  
 \item $f$ is an $\FF$-definable{\footnote{That is, $f$ is definable in the real ordered field $(\RR,+,\cdot, 0,1,\leq)$ by a first-order formula with parameters from $\FF$.}}, analytic function that is one-to-one in each coordinate;
 
 \item $f(t_0,t_1, \ldots, t_{d-1}) = a$.
 
\end{itemize}
For each $a \in \RR^n$, let $f_a$ be a determining function for $a$. 
Since each $f_a$ is $\FF$-definable and $|\FF| = \kappa$, then $|\{f_a : a \in \RR^n\}| \leq \kappa$.
For each $a \in \RR^n$, let $A_a = \{b \in \RR^n : f_a$ is a determining function for $b\}$. 
Let ${\mathcal A} = \{A_a : a \in \RR^n\}$. Clearly, $|{\mathcal A}| \leq \kappa$ and $\RR^n = \bigcup {\mathcal A}$. 

Consider some $A_a \in {\mathcal A}$. Let $d < \omega$ and 
$\dom(f_a) = (q_0, r_0) \times (q_1, r_1) \times \cdots \times   (q_{d-1}, r_{d-1})$,
For $i < d$, let $X_i = (q_i,r_i) \cap T$. Then $|X_0| = |X_1| = \cdots = |X_{d-1}| = 2^{\aleph_0}$.
Let $X = X_0 \times X_1 \times \cdots \times X_{d-1}$. Then $f_a \harp X$ is a bijection from 
$X$ to $A_a$. For $i < d$, let $h_i : X_i \into \RR$ be a bijection, and then define $h : X \into \RR^n$
so that $h(x_0,x_1, \ldots, x_{d-1}) = \langle h_0(x_0), h_1(x_1), \ldots, h_{d-1}(x_{d-1}) \rangle$. Hence, 
 $h : X \in \RR^d$ is a bijection. Now let $g = hf_a^{-1}$; hence,  $g : A_a \into \RR^d$ is a bijection.

Now consider some $H \in {\mathcal H}$. Let $F \subseteq [X]^k$ be such that $f_a$ is an isomorphism from $(X,F)$ to $H(A_a)$. Each edge of $P \in F$ is a $d$-dimensional $k$-template, and   every homomorphic image $Q \in [X]^k$ is also in $F$. Thus, $(X,F)$ is a multi-template hypergraph, so there is a set ${\mathcal P}$ of $d$-dimensional $k$-templates such that $(X,F) = L(X,{\mathcal P})$.
Then, $h$ is an isomorphism from $L(X,{\mathcal P})$ to $(\RR^d, {\mathcal P})$, so $g$ is an isomorphism from $H[A_a]$ to $L(\RR^d, {\mathcal P})$. \qed 
 
 \bigskip
%
%

{\bf \S2.~The chromatic numbers of algebraic hypergraphs.}
If $H = (V,E)$ is a $k$-hypergraph, then a function $\varphi : V \into C$ is a \mbox{$C$-{\bf coloring}} of $H$. The $C$-coloring $\varphi$ is {\bf proper} if whenever $\{a_0,a_1, \ldots,$ $ a_{k-1}\} \in E$, then there are $i < j < k$ such that $\varphi(a_i) \neq \varphi(a_j)$. The chromatic number $\chi(H)$ of $H$ is the least (finite or infinite) cardinal $\kappa$ such that there is a proper $\kappa$-coloring of $H$. 
 The purpose of \cite{sch16} was to determine the chromatic numbers of all algebraic hypergraphs 
 with infinite chromatic numbers.  This improved \cite{sch00} in which  those algebraic hypergraphs having uncountable chromatic numbers were characterized. 
 
 If $\kappa$ is an infinite cardinal, then $\kappa^+$ is its successor cardinal. For $n < \omega$, we define 
 $\kappa^{+n}$ by recursion so that $\kappa^{+0} = \kappa$ and $\kappa^{+(n+1)} = (\kappa^{+n})^+$.
 We define $\kappa^{-n}$ to be the least cardinal $\lambda$ such that $\lambda^{+n} \geq \kappa$. 
 We make use of the symbol $\infty$: for any infinite $\kappa$,  $\kappa^{+\infty} = \infty > \kappa$ and $\kappa^{-\infty} = \aleph_0$.  
 
 The chromatic numbers of some template hypergraphs are given in the next lemma. If $P$ is a $d$-dimensional $k$-template, then a subset $I \subseteq d$ is a {\bf distinguisher} for $P$ if 
  whenever $x,y \in P$ are distinct, then $x_i \neq y_i$ for some $i \in I$. We let $e(P)$ be the least cardinality of a distinguisher for $P$.

  \bigskip
  
  {\sc Lemma} 2.1: (\cite[Th.~2.1]{sch16}) {\em If $P$ is a $d$-dimensional $k$-template and $X$ is an infinite set, then
  $$
  \chi(L(X^d,P)) = |X|^{-(e(P) - 1)}. 
  $$}
  
  \bigskip
  
  {\sc Corollary} 2.2: {\em Suppose that ${\mathcal P}$ is a set of  $d$-dimensional $k$-templates and $X$ is an infinite set. If ${\mathcal P} = \varnothing$, then $\chi(L(X^d,{\mathcal P})) = 1$. If ${\mathcal P} \neq \varnothing$, then 
  $
  \chi(L(X^d,{\mathcal P})) = |X|^{-(e - 1)}
  $, where $e = \min\{e(P) : P \in {\mathcal P}\}$.} \qed
  
  \bigskip

 There is a certain function $\delta$ that we will use such that $\delta(p(\overline x)) \in \omega \cup \{\infty\}$ for every $(k,n)$-ary polynomial $p(\overline x)$. We leave this function undefined for now, but do note that, by results in \cite{sch00}, 
 it is (much more than) absolute. If $p(\overline x)$ defines $H$ and $\chi(H)$ is finite, then $\delta(p(\overline x)) = \infty$. Thus, Theorem~2.4 could be used as a definition. It is possible that $\chi(H)$ is infinite and $\delta(p(\overline x)) = \infty$. For example,
 $\delta(p_n(x,y,z,w)) = \delta(q_n(x,y,z)) = \delta(p_n(x,y,z,w) \cdot q_n(x,y,z)) = \infty$ even though the hypergraphs they define have chromatic number $\aleph_0$.  It is proved in \cite{sch17} that when restricted to polynomials with rational coefficients, $\delta$ is computable. In fact, in a sense that is given in \cite{sch17}, the function $\delta$ is decidable.
 
 \bigskip

  The  following proposition, which  is relevant
to (1) and (2) of Theorem~2, says exactly when $\delta(p(\overline x)) = 0$. 
  
  \bigskip
  
  {\sc Proposition 2.3}: {\em If $p(\overline x)$ is a $(k \times n)$-ary polynomial that defines the hypergraph 
   $H$, then $\delta(p(\overline x)) = 0$ iff $H$ has an infinite complete subset.}
   
   \bigskip
 
 Thus, an additional equivalent could be added to Theorem~2: if $p(\overline x)$ defines $H$, then 
 $\delta(p(\overline x)) > 0$.  The following is the main result of \cite{sch16}. 
 
 \bigskip
 
 {\sc Theorem 2.4}: {\em If $p(\overline x)$ defines $H$ and $\chi(H)$ is infinite,  then, $\chi(H) = (2^{\aleph_0})^{-\delta(p(\overline x))}$.}
 
 \bigskip
 
 A consequence of Lemma~2.1 and Theorem~2.4 is: If $p(\overline x)$ is a polynomial that defines 
 the template $k$-hypergraph $L(\RR^d,P)$, then $e(P) = 1 +\delta(p(\overline x))$.

 One half of the conclusion of Theorem~2.4 is the existence of a proper $(2^{\aleph_0})^{-\delta(p(\overline x))}$-coloring of $H$.  The proof of Theorem~1.2 in \cite{sch16} shows that often a single coloring (the {\em master coloring}) can do the job of many. The following theorem is not explicitly stated in \cite{sch16}, 
 but the  special case (\cite[Coro.~3.7]{sch16} of it that concerns $D$-distance graphs is. 
 
 \bigskip
 
 {\sc Theorem 2.5}: {\em Suppose $1 \leq n < \omega$,  $\aleph_0 \leq \kappa \leq 2^{\aleph_0}$ and ${\mathcal H}$ is a set 
 of algebraic hypergraphs $H = (\RR^n,E)$ such that $|{\mathcal H}| \leq \kappa$ and  $\chi(H) \leq \kappa$ for each $H \in {\mathcal H}$. Then there is $\varphi : \RR^n \into \kappa$ that is a proper $\kappa$-coloring of each $H \in {\mathcal H}$. }
 
 \bigskip
 
 {\it Proof}. 
 Once we have Theorem~1.2, the proof is  straightforward.  Let ${\mathcal A}$ be as in Theorem~1.2. We will define $\varphi : \RR^n \into C$, where $|C| = \kappa$, that is a proper coloring of each $H \in {\mathcal H}$. Let $a \in \RR^n$. Choose some $A \in {\mathcal A}$ such that $a \in A$.  Then we have $d < \omega$ and $g : A \into \RR^d$ as in Theorem~1.2. Let $\theta : \RR^n \into \kappa$ be such that it is a proper coloring of each of the finitely many multi-template hypergraphs $L(\RR^d, {\mathcal P})$ for which $\chi(L(\RR^d, {\mathcal P})) \leq \kappa$. Then let 
 $\varphi(a) = \langle A, \theta g(a) \rangle$. One easily checks that $\varphi$ is as intended. \qed

 \bigskip
 
{\sc Theorem 2.6}: {\em If $p(\overline x)$ defines $H = (\RR^n,E)$ and $X \subseteq \RR^n$ is infinite. Then $\chi(H[X]) \leq  |X|^{-\delta(p(\overline x))}$.} 

\bigskip

{\it Proof}. Let $H$ be a $k$-hypergraph. We can assume that $\chi(H)$ is infinite, as otherwise $\chi(H[X])$ is also finite. Hence,  $\chi(H) = (2^{\aleph_0})^{-\delta(p\overline x))}$  by Theorem~2.4.  Let ${\mathcal A}$ be as in Theorem~1.1. Thus, for  each $A \in {\mathcal A}$, $H[A]$ is isomorphic to some multi-template $k$-hypergraph. 
By Corollary~2.2, for each $A$, either $\chi(H[A]) = 1$ or $\chi(H[A]) = (2^{\aleph_0})^{-(e_A-1)}$, where $1 \leq e_A < \omega$.  Thus, each $e_A \geq \delta(p(\overline x))$. 
 It then follows from Corollary~2.2 again, that for any infinite $Z \subseteq A$,
$\chi(H[Z]) \leq |Z|^{-(e_A-1)}$. Thus, $\chi(H(X \cap A)) \leq |X|^{-(e_A-1)} \leq |X|^{-\delta(p(\overline x))}$. Since ${\mathcal A}$ is countable, it follows that $\chi(H[X]) \leq  |X|^{-\delta(p(\overline x))}$. \qed

\bigskip

 As a postscript,   there is the following proposition that asserts that the given bound in Theorem~2.6  is optimal.
 
 \bigskip
 
 {\sc Proposition 2.7}: {\em  If $p(\overline x)$ defines $H = (\RR^n,E)$, $\chi(H) \geq \aleph_0$ and $\aleph_0 \leq \kappa \leq 2^{\aleph_0}$, then there is $X \subseteq \RR^n$  such that $|X| = \kappa$ 
 and  $\chi(H[X]) =  \kappa^{-\delta(p(\overline x))}$.} 

\bigskip

{\it Proof.} Let $p(\overline x), H$ and $\kappa$ be as given. Let ${\mathcal A}$ be as in Theorem~1.1. 
Let $A \in {\mathcal A}$ be such that $\chi(H) = \chi(H[A])$. There are an algebraic multi-template 
hypergraph $L(\RR^d, {\mathcal P})$ such that $\chi(H) = \chi(L(\RR^d, {\mathcal P}))$ and an isomorphism $g$ from $H[A]$ to $L(\RR^d, {\mathcal P})$. Since $\chi(H)$ is infinite, ${\mathcal P} \neq \varnothing$. Let $e = \min\{e(P) : P \in {\mathcal P}\}$. Then $e = \delta(p(\overline x)) + 1$. 
Let $Y \subseteq \RR$ such that $|Y| = \kappa$. Then $\chi(L(Y^d,{\mathcal P}) =  \kappa^{-\delta(p(\overline x))}$. Let $X = g^{-1}(Y^d)$. Then $|X| = \kappa$ and $\chi(H[X]) = \kappa^{-\delta(p(\overline x))}$. \qed

\bigskip

{\sc Theorem 2.8}: {\em Suppose $1 \leq n < \omega$,  $\aleph_0 \leq \kappa \leq 2^{\aleph_0}$, $X \subseteq \RR^n$, $|X| = \kappa$ and ${\mathcal H}$ is a set 
 of algebraic hypergraphs $H = (\RR^n,E)$ such that $|{\mathcal H}| \leq \kappa$ and  $\chi(H[X]) \leq \kappa$. Then there is $\varphi : \RR^n \into \kappa$ that is a proper $\kappa$-coloring of each $H \in {\mathcal H}$.}
 
 \bigskip
 
 {\it Proof}. This proof is just like the proof of Theorem~2.6, but using Theorem~1.2 instead of Theorem~1.1. \qed
 
 \bigskip
 
%
%

{\bf \S3.~The proof of Theorem~2}.  This section is devoted to completing the proof of Theorem~2.

Trivially, $(1) \Longrightarrow (2)$. 
The converse implication $(2) \Longrightarrow (1)$ is a consequence of \cite[Coro.\@ 4.9]{sch17}, 
in which $(2)$ is replaced by the seemingly weaker

\begin{itemize}

\item[($2'$)] {\em If $X \subseteq \RR^n$ is algebraic and complete for $p(x_0,x_1, \ldots,x_{k-1})$, then  $X$ is countable.}

\end{itemize} 
It is obvious that $(3) \Longrightarrow (2)$. For, if (2) fails for some $X$, then this $X$ is also a counter-example to $(3)$. 
Thus, it remains to prove $(1) \Longrightarrow (3)$. 
We will prove  an apparently stronger theorem which involves  many hypergraphs simultaneously rather than just one. 

\bigskip

{\sc Theorem 3.1}: {\em Suppose that $1 \leq n < \omega$ and $\aleph_0 \leq \kappa  \leq 2^{\aleph_0}$.
Suppose that ${\mathcal H}$ is a set of  algebraic hypergraphs $(\RR^n,E)$ such that $|{\mathcal H}| < \kappa$ and 
no $H \in {\mathcal H}$ has an infinite complete subset. Then, every $X \subseteq \RR^n$  such that $|X| = \kappa$ has a subset $Y \subseteq X$ such that $|Y| = |X|$ and $Y$ is independent for every $H \in {\mathcal H}$.}

\bigskip

{\it Proof}. Fix $1 \leq n < \omega$.  We will prove the theorem by induction on~$\kappa$. 

\smallskip

$\kappa = \aleph_0$: Let ${\mathcal H}$ and $X$ be as given. By several applications of Ramsey's Theorem, one for each $H \in {\mathcal H}$, there is an infinite $Y \subseteq X$ such that for each $H \in {\mathcal H}$, either $Y$   is  independent for $H$ or $Y$ is complete for $H$. Since no $H \in {\mathcal H}$ has an infinite complete subset,  $Y$ is independent for each $H \in {\mathcal H}$.
\smallskip

$\kappa > \aleph_0$: Suppose that for all smaller cardinals, the theorem is true. Let ${\mathcal H}$ and $X$ be as given. 
There are two cases depending on whether or not $\kappa$ is a successor cardinal.

\smallskip

{\em $\kappa$ is a successor cardinal}: Let $\kappa = \lambda^+$. For each $H \in {\mathcal H}$, let $p_H(\overline x)$ be a polynomial that defines $H$ and let $e_H = e(p_H(\overline x))$.  Since no $H \in {\mathcal H}$ has an infinite complete set, then by Proposition 2.3, each $e_H \geq 1$. Thus,  by Theorem~2.6, $\chi(H[X]) \leq \kappa^{-\delta(p_H(\overline x))} \leq \lambda$.  By Theorem~2.8, there is $\varphi : X \into \lambda$ that is a proper coloring of $H[X]$ for each $H \in {\mathcal H}$. Since $\kappa$ is regular, we let $Y \subseteq X$ be such that $|Y| = \kappa$ and $\varphi$ is constant on $Y$. Then, $Y$ is independent for every $H \in {\mathcal H}$. 

\smallskip

{\em $\kappa$ is a limit cardinal}: Let $\lambda = \cf(\kappa)$, and let $\langle \kappa_\alpha : \alpha < \lambda \rangle$ be a sequence of regular cardinals and $\langle X_\alpha : \alpha < \lambda \rangle$  a sequence of subsets of $X$ such that:

\begin{itemize}

\item if $\alpha < \lambda$, then $|X_\alpha| = \kappa_\alpha < \kappa$;

\item $\kappa = \bigcup_{\alpha < \lambda} \kappa_\alpha$;

\item if $\alpha  < \lambda$, then $\kappa_\alpha > \bigcup_{\beta< \alpha} \kappa_\beta$;

\item if $\alpha < \beta < \lambda$, then $X_\alpha \cap X_\beta = \varnothing$;

\item  $\kappa_0 > |{\mathcal H}|$.

\end{itemize}

For each $H \in {\mathcal H}$, let $p_H(x_0,x_1, \ldots, x_{k_H-1})$ be a polynomial that defines $H$. 
For each $\alpha < \lambda$, let ${\mathcal H}_\alpha$ be the set of all hypergraphs defined by 
a polynomial obtained from some $p_H(x_0,x_1, \ldots, x_{k_H-1})$ by substituting some elements of 
$\bigcup_{\beta < \alpha}X_\beta$ for some (but not all) of the variables $x_0,x_1, \ldots, x_{k_H-1}$.
For example, ${\mathcal H}_0 = {\mathcal H}$, so $|{\mathcal H}_0| < \kappa_0$.  If  $0 < \alpha < \lambda$, then 
$|{\mathcal H}_\alpha| \leq \bigcup_{\beta < \alpha}\kappa_\beta $. Thus, $|{\mathcal H}_\alpha| < \kappa_\alpha$ for all $\alpha < \lambda$.
Notice that no hypergraph in any ${\mathcal H}_\alpha$ has an infinite complete subset. 

By the inductive hypothesis, we can get $Y_\alpha \subseteq X_\alpha$, for each $\alpha < \lambda$, 
such that $|Y_\alpha| = |X_\alpha|$ and $Y_\alpha$ is independent for all hypergraphs in ${\mathcal H}_\alpha$. Let $Y = \bigcup_{\alpha < \lambda}Y_\alpha$. Then $Y \subseteq X$ and $|Y| = |X|$. We claim that 
$Y$ is independent for each $H \in {\mathcal H}$. The claim is easily proved by showing by induction on $\alpha < \lambda$ that  $\bigcup_{\beta < \alpha}Y_\beta$ is independent for each hypergraph in ${\mathcal H}_\alpha$. 
 \qed

\bigskip 

We conclude with a sample corollary of Theorem~3.1.

\bigskip

{\sc Corollary 3.2}: {\em Suppose that $1 \leq n < \omega$. Every uncountable $X \subseteq \RR^n$ has a 
subset $Y \subseteq X$ such that $|Y| = |X|$ and $Y$ has no algebraically dependent distances.} \qed

\bigskip

To clarify the conclusion of this corollary, it means: whenever $\{a_0,b_0\},$ $ \{a_1,b_1\}, \ldots \{a_m,b_m\} \in [Y]^2$ are distinct, then   
$\|a_m-b_m\|$ is not algebraic over $\{\|a_i-b_i \| : i < m\}$.

\bibliographystyle{plain}

\end{document}